\documentclass[letterpaper, 10pt, journal]{IEEEtran}
\usepackage{cite}
\usepackage{amsmath} 
\usepackage{amssymb}  
\usepackage{amsfonts}
\usepackage{mathrsfs}
\usepackage{mathtools}
\usepackage[dvipsnames]{xcolor}
\usepackage[thmmarks, amsmath]{ntheorem}
\usepackage{enumerate}
\usepackage{booktabs}
\usepackage[bbgreekl]{mathbbol}
\usepackage{colortbl}
\usepackage{ragged2e}
\usepackage{enumitem}
\usepackage{relsize}
\usepackage{standalone}
\usepackage{booktabs, array}
\usepackage{siunitx}
\usepackage{url}
\usepackage[ruled,vlined]{algorithm2e}
\graphicspath{{./data/}}

\usepackage[font=footnotesize]{caption}

\usepackage{graphicx}

\newcommand{\norm}[1]{\ensuremath{\| #1 \|}}
\newcommand{\until}[1]{[#1]}

\newcommand{\setdef}[2]{\{#1 \; | \; #2\}}

\newcommand{\MM}{\mathcal{M}}

\newcommand{\cl}{\operatorname{cl}}

\newcommand{\st}{\operatorname{s.} \operatorname{t.}}

\newcommand{\Eb}{\mathbb{E}}
\newcommand{\Pb}{\mathbb{P}}
\newcommand{\Rb}{\mathbb{R}}

\usepackage{subcaption}

\newcommand{\GG}{\mathcal{G}}

\newcommand{\RR}{\mathcal{R}}
\newcommand{\LL}{\mathcal{L}}

\newcommand{\TT}{\mathcal{T}}

\newcommand{\PP}{\mathcal{P}}

\newcommand{\CVaR}[1]{\operatorname{CVaR}^{#1}}

\newcommand{\Pbhat}{\widehat{\Pb}}
\newcommand{\data}{\widehat{\omega}}

\usepackage{siunitx}

\newcommand{\ignore}[1]{}
\allowdisplaybreaks
\DeclareSymbolFont{bbold}{U}{bbold}{m}{n}
\DeclareSymbolFontAlphabet{\mathbbold}{bbold}

\newcommand{\oprocendsymbol}{\hbox{$\bullet$}}
\newcommand{\oprocend}{\relax\ifmmode\else\unskip\hfill\fi\oprocendsymbol}

\newtheorem{theorem}{Theorem}[section]

\newtheorem{proposition}[theorem]{Proposition}

\newtheorem{remark}[theorem]{Remark}

\newcolumntype{R}[1]{>{\RaggedLeft\arraybackslash}p{#1}}

\newcolumntype{P}[1]{>{\centering\arraybackslash}p{#1}}

\makeatletter
\def\subsection{\@startsection{subsection}{2}{\z@}{1ex plus 1ex minus 0ex}%
{0.7ex plus .5ex minus 0ex}{\normalfont\normalsize\itshape}}%
\makeatother

\usepackage[absolute,overlay]{textpos}
\usepackage[hyperindex=true, %
  colorlinks=true,%
  pagebackref=false,%
  plainpages=false,%
  pdfpagelabels,%
  linkcolor=NavyBlue,%
  citecolor=NavyBlue,%
  filecolor=NavyBlue,%
  urlcolor=NavyBlue%
]{hyperref} 

\hypersetup{
	pdftitle={Wasserstein Distributionally Robust Look-Ahead Economic Dispatch},
	pdfauthor={Bala Kameshwar Poolla, Ashish R. Hota, Saverio Bolognani, Duncan S. Callaway, and Ashish Cherukuri}
	}

\begin{document}

\begin{textblock*}{\textwidth}(15mm,9mm) % {block width} (coords) 
\centering \bf \textcolor{NavyBlue}{To appear in the \emph{IEEE Transactions on Power Systems} \\\url{https://doi.org/10.1109/TPWRS.2020.3034488}}
\end{textblock*}

\title{Wasserstein Distributionally Robust \\ Look-Ahead Economic Dispatch}
\author{Bala Kameshwar Poolla, Ashish R. Hota, Saverio Bolognani, Duncan S. Callaway, and Ashish Cherukuri
\thanks{Bala Kameshwar Poolla and Duncan S. Callaway are with the Energy and Resources Group, University of California Berkeley, CA. Email: {\tt \{bpoolla, dcal\}@berkeley.edu}}%
\thanks{Ashish R. Hota is with the Indian Institute of Technology Kharagpur, India. Email: {\tt ahota@ee.iitkgp.ac.in}}%
\thanks{Saverio Bolognani is with the Automatic Control Laboratory, ETH Zurich, 8092 Zurich, Switzerland. Email: {\tt bsaverio@ethz.ch}}%
\thanks{Ashish Cherukuri is with the University of Groningen, Netherlands. Email: {\tt a.k.cherukuri@rug.nl}}%
\thanks{This work was partially supported by a grant from IIT Kharagpur under the ISIRD scheme and by the NSF under the Cyber SEES-1539585.}%
\thanks{\textcopyright 2020 IEEE.  Personal use of this material is permitted.  Permission from IEEE must be obtained for all other uses, in any current or future media, including reprinting/republishing this material for advertising or promotional purposes, creating new collective works, for resale or redistribution to servers or lists, or reuse of any copyrighted component of this work in other works.}%
}

\thispagestyle{plain}
\pagestyle{plain}
\maketitle

\begin{abstract}
We consider the problem of look-ahead economic dispatch (LAED) with uncertain renewable energy generation.
The goal of this problem is to minimize the cost of conventional energy generation subject to uncertain operational constraints.
The risk of violating these constraints must be below a given threshold for a family of probability distributions with characteristics similar to observed past data or predictions.
We present two data-driven approaches based on two novel mathematical reformulations of this distributionally robust decision problem. The first one is a tractable convex program in which the uncertain constraints are defined via the distributionally robust conditional-value-at-risk. The second one is a scalable robust optimization program that yields an approximate distributionally robust chance-constrained LAED. Numerical experiments on the IEEE 39-bus system with real solar production data and forecasts illustrate the effectiveness of these approaches. We discuss how system operators should tune these techniques in order to seek the desired robustness-performance trade-off and we compare their computational scalability.
\end{abstract}

\section{Introduction}
\label{section:introduction}
The electricity grid is witnessing an increasing penetration of renewable energy sources (such as solar photovoltaic and wind) \cite{kroposki2017achieving}. In sharp contrast with conventional sources of electricity (such as coal-fired or nuclear power plants), the energy produced from renewable energy sources is highly variable, intermittent, and not fully dispatchable. Thus, efficient integration of renewable energy sources so as to meet the demand for electricity while respecting the operational constraints (such as line flow limits and ramp constraints) is a fundamental challenge for modern power grids \cite{bienstock2014chance}.\\
\indent The problem of determining the cost-efficient dispatch schedule for (conventional) generators in order to meet the forecast demand subject to operational constraints is referred to as the {\it optimal power flow (OPF)} or {\it economic dispatch (ED)} problem \cite{bienstock2014chance,varaiya2010smart}. Both single-period as well as multi-period version $-$ referred to as the {\it look-ahead economic dispatch (LAED)}, have been investigated \cite{ross1980dynamic,varaiya2010smart}.\\ 
\indent {In this work, we investigate the multi-period LAED problem in the presence of uncertain renewable energy generation. In practical terms, the multi-period decision process allows the operator to take strategic decisions (such as ramping-up the most economical conventional generation) several hours ahead of real-time operation in order to ensure sufficient controllability of the system for a range of possible realizations of the uncertain renewable generation \cite{capitanescu2011state}.}\\
\indent This dispatch problem under uncertain renewable energy generation results in a robust or a stochastic optimization problem. The decision-maker either requires the uncertain constraints to hold for all realizations of the uncertainty (leading to a robust/worst-case optimization formulation) or with a high probability (leading to a chance-constrained program) \cite{shapiro2009lectures}. The latter yields less conservative solutions, but requires the decision-maker to know the distribution of the uncertainty. {Some early works have modeled the distribution of (forecast error in) wind power generation as Gaussian \cite{albadi2011comparative} or Beta \cite{bludszuweit2008statistical}. However, both these models were challenged in subsequent works, e.g., \cite{tewari2011statistical}. Other distributions, such as Laplace \cite{wu2014statistical}, Cauchy \cite{hodge2011wind}, and Levy alpha-stable \cite{bruninx2014statistical} were also proposed. Nevertheless, as discussed in \cite{wang2015distributionally}, there is no probability distribution that is suitable to describe all wind energy generation data. Analogous observations have been made in \cite{golestaneh2016very} regarding solar energy generation. Furthermore, climate change also induces subtle shifts in renewable energy generation compared to historical data \cite{solaun2019climate}.}\\ 
\indent {This lack of a suitable distribution that describes renewable energy generation on one hand, and availability of historical and numerical forecast data on the other, have been one of the primary motivations behind the rise of distributionally robust approaches to solve various operational problems in modern power systems. In a distributionally robust chance-constrained program (DRCCP), the goal is to find a solution which satisfies the chance-constraints for a suitably defined family of distributions of the uncertain parameters (as opposed to learning a single distribution that best captures the observed data and requiring the chance-constraint to be satisfied for this learned distribution). The family of distributions is referred to as an {\it ambiguity set}. Thus, this approach enables the decision-maker to robustify dispatch decisions to slow trends, seasonal variations, and non-ergodicity in the renewable generation data and avoid overfitting to observed data.}

\subsection{Related Works} 

Early work on distributionally robust (DR) OPF considered moment-based ambiguity sets which comprise of all distributions with an identical mean and covariance as the uncertain parameters \cite{summers2015stochastic,zhang2016distributionally,li2019distributionally}. However, this requires to infer the mean and covariance of the uncertain parameters from the empirical data in order to construct the ambiguity set. In contrast, recent papers have considered the Wasserstein distributionally robust optimization paradigm in power systems applications such as unit commitment\footnote{{The unit commitment problem is an instance of an integer program which belongs to a different class of optimization problems than the OPF or LAED problem considered here (with continuous decision variables).}} \cite{duan2017data,zhu2019wasserstein} and optimal power flow \cite{guo2018data,wanggao2018,duan2018distributionally}.

The definition of ambiguity sets via the Wasserstein distance, directly utilizing the observed samples, brings several attractive properties in terms of finite sample guarantees, tractable reformulations, and asymptotic consistency.
{In Section~\ref{ssec:discussion}, we briefly review the underlying assumptions and rigorous finite sample guarantees established for Wasserstein ambiguity sets in prior work \cite{fournier2015rate,weed2019sharp,esfahani2018data} (and also for the scenario approach, which requires milder assumptions).} 

{Our contribution builds upon earlier works, in particular \cite{esfahani2018data,guo2018data,duan2018distributionally}. The authors in \cite{esfahani2018data} were the first to propose finite-dimensional reformulations of distributionally robust optimization problems with uncertain cost functions over Wasserstein ambiguity sets; however \cite{esfahani2018data} does not deal with distributionally robust chance or risk-constrained optimization problems. In \cite{guo2018data}, the reformulations developed in \cite{esfahani2018data} were applied to the multi-period OPF problem. While \cite{guo2018data} (as well as \cite{wanggao2018}) notes that the operational constraints (such as line flow and voltage magnitude limits) in the OPF problem are uncertain under renewable energy generation, they treat these constraints as penalty terms in the cost function. This allows them to use the results of \cite{esfahani2018data}. However, handling constraints by moving them to the objective function via penalty terms does not guarantee that the constraints will be satisfied at the optimum. 
In \cite{duan2018distributionally}, uncertain operational constraints are treated as distributionally robust individual chance-constraints under Wasserstein ambiguity sets. As a consequence of their modeling choice, on each constraint, the uncertainty takes the form of a scalar random variable, which is the basis for their reformulations. However, solutions obtained under individual chance-constraints lack the desired robustness of the solutions obtained under joint chance-constraints.}

\subsection{Summary of Contributions} 
{In this work, we study the LAED problem where the operational constraints are required to hold jointly with a high probability for all distributions that are ``close" to the empirical distribution induced by the observed data or by the available forecasts, as measured by the Wasserstein metric. However, chance-constrained programs are computationally intractable except for a special class of distributions and constraints, even when the distribution of the uncertain parameters is known.} Accordingly, past work has focused on developing tractable convex approximations of the chance-constrained OPF problem \cite{summers2015stochastic,baker2019joint,halilbavsic2018convex}. Sample based methods, inspired by the so-called {\it scenario approach} \cite{campi2008exact} and its variations, have also been investigated in this context \cite{vrakopoulou2013probabilistic,modarresi2018scenario,chamanbaz2019probabilistically}. We adopt a similar approach here and develop two tractable approximations for Wasserstein DRCCPs.

First, we observe that conditional-value-at-risk (CVaR)-constraints act as convex inner approximations to chance-constraints \cite{nemirovski2006convex}. Furthermore, CVaR is a widely used coherent risk measure \cite{rockafellar2000optimization} which guarantees that the constraints not only hold with high probability, but also the magnitude of constraint violation is small in expectation. {We present a convex finite-dimensional reformulation of distributionally robust CVaR-constrained programs (DRCVP) under Wasserstein ambiguity sets for constraint functions that are affine in the decision variables and the uncertain parameters.\footnote{{This reformulation appeared in a preliminary version of this work \cite{hota2019data} without proof. Here, we include the complete proof of this result. We show that the duality results derived in \cite{esfahani2018data} are not directly applicable for our problem, but under relatively mild conditions, a finite-dimensional tractable reformulation can be derived.}}}

However, the number of constraints of the DRCVP problem increases with the number of samples, leading to high dimensionality, despite the convexity. {Therefore, we develop a scalable approach to approximately solve DRCCPs} under Wasserstein ambiguity sets inspired by a similar approach proposed in \cite{margellos2014road} for chance-constrained programs. {We leverage a recently proposed exact reformulation of DRCCPs in \cite{hota2019data} 
and approximate the problem via a two-dimensional DRCCP for each component of the uncertainty and a master robust optimization problem whose size does not depend on the number of samples.

These tractable reformulations apply not only to the LAED problem, but to any distributionally robust chance or CVaR-constrained optimization problem over Wasserstein ambiguity sets, as long as the constraint function is affine in the decision variables and uncertainty. We present rigorous proofs of our theoretical results which make the paper self-contained.

Finally, we carry out an extensive empirical evaluation of the proposed formulations by solving the LAED problem for the IEEE $39$-bus transmission grid with real solar irradiation and forecast data. In particular, we consider two settings:
\begin{itemize}[leftmargin=0.2cm, itemindent=0.5cm]
\setlength{\itemsep}{1pt}
\item the operator has access to an ensemble of forecasts of the solar irradiation for the next day;
\item the operator has access to past data on solar generation.
\end{itemize}
{For both settings, we discuss why the theoretical guarantees available in the literature cannot be used to tune the size of the ambiguity set, specifically because the renewable generation data does not consist of i.i.d. samples from an underlying distribution (see Section \ref{sec:fsg} for a detailed discussion). Instead, we empirically evaluate the trade-off between robustness (out-of-sample constraint satisfaction) and performance as a function of the size of the ambiguity set.\footnote{{Prior applications of Wasserstein distributionally robust optimization in power systems such as \cite{guo2018data2,zhu2019wasserstein} have also primarily relied on empirical evaluation of the trade-off between robustness and performance.}}} Finally, we discuss how these approaches scale with network size and highlight several interesting avenues for future research.
\section{Look-Ahead Economic Dispatch under Uncertainty}
\label{section:problem_formulation}

In this section, we define the look-ahead economic dispatch (LAED) problem. Our formulation is inspired by a similar structure in \cite{modarresi2018scenario}. The objective of the LAED problem is to minimize the total cost of generation over a time-horizon of length $T$, while satisfying operational constraints and the forecasted power demand in an appropriate sense in the presence of uncertain power generation from renewable (solar and wind) energy sources (RESs). With a slight abuse of notation, let $\GG, \RR$, and $\LL$ denote the set of conventional generators, RESs, and loads as well as the corresponding nodes in the power network, and let $|\GG| = N_g, |\RR| = N_r$, and $|\LL| = N_\ell$. The sets $\GG$ and $\RR$ need not be disjoint, i.e., a node may have both conventional generation and RESs. 

We denote by $p_{i}[t], w_{j}[t]$, and $\ell_{k}[t]$ the power generation of the conventional generator $i \in \GG$, the RES $j \in \RR$, and the power consumed by a load $k \in \LL$ at time $t$, respectively. The corresponding aggregate quantities in vector form are denoted by $p[t] \in \Rb^{N_g}$, $w[t] \in \Rb^{N_r}$, and $\ell[t] \in \Rb^{N_\ell}$, respectively. 
Now, let $c_{i}[t]$ denote the per-unit cost of power generation for the conventional power plant $i \in \GG$ at time $t$. We assume that the marginal cost of renewable energy generation is zero. {Let $\TT := \{t_0+1,\ldots,t_0+T\}$ be the optimization horizon (for example, the day ahead setting). 
We assume that the generation set-points $p[t_0]$ at the starting time $t_0$ are known (as they are part of today's schedule).
The LAED problem with starting time $t_0$ and horizon $T$ is mathematically expressed as%
\begin{subequations}\label{eq:LAED}
\begin{align}
\min_{\{p[t]\}_{t \in \TT}} & \, \, \sum_{t \in \TT}\,\sum_{i \in \GG} c_{i}[t] p_{i}[t]\\
\st & \, \, \mathrm{RD}_i[t] \leq p_{i}[t]-p_{i}[t-1] \leq \mathrm{RU}_i[t], \forall i \in \GG,  \label{eq:ramp_constraints} \\
&  \, \, \underline{P_{i}}[t] \leq p_{i}[t] \leq \overline{P_{i}}[t], \forall i \in \GG, \label{eq:power_bound} \\
&  \, \, \sum_{i \in \GG} p_{i}[t] + \sum_{j \in \RR} w_{j}[t] \geq \sum_{k\in \LL} \ell_{k}[t], \label{eq:powerbalance_stochastic} \\
&  \, \, -\overline{F} \! \leq \! \Lambda \Big[B_p p[t] + B_w w[t] - B_\ell \ell[t]\Big] \! \leq \! \overline{F}, \label{eq:lineflow_constraint}
\end{align}
\end{subequations}
where the constraints \eqref{eq:ramp_constraints}, \eqref{eq:power_bound}, \eqref{eq:powerbalance_stochastic}, and \eqref{eq:lineflow_constraint} hold for $t \in \TT$.}

The power generation of the conventional generators are the decision variables, the power consumed by the loads are assumed to be known, and the power generation of the RESs are treated as uncertain parameters. The parameter $\mathrm{RD}_i[t]$ (respectively, $\mathrm{RU}_i[t]$) denotes the ramp-down (respectively, ramp-up) capacity and $\underline{P_{i}}[t]$ (respectively, $\overline{P_{i}}[t]$) denotes the lower bound (respectively, upper bound) of the conventional generator $i$ at time $t$. Thus, the constraints \eqref{eq:ramp_constraints} and \eqref{eq:power_bound} are deterministic. Although written as deterministic for ease of representation, \eqref{eq:powerbalance_stochastic} and~\eqref{eq:lineflow_constraint} have uncertain or stochastic parameters $w[t], t \in \TT$. The constraint \eqref{eq:powerbalance_stochastic} ensures sufficient generation and can be adapted to account for the available reserves ({see Remark \ref{remark:powerbalance})}.

The inequality \eqref{eq:lineflow_constraint} requires line power flows to be within permissible limits with the vector of line flow limits denoted by $\overline{F}$. The flows in the transmission lines are computed by leveraging the so-called Power Transfer Distribution Factor (PTDF) matrix $\Lambda$, a linear sensitivity that represents the marginal change of the active power flow on a line  if we apply a marginal increase of the power injection at a node. 
More specifically, let $N_s$ and $N_e$ be the total number of nodes and lines (edges) in the network and $B_g \in \Rb^{N_s\times N_g}$, $B_r \in \Rb^{N_s\times N_r}$, and $B_\ell \in \Rb^{N_s\times N_\ell}$ denote appropriate matrices.\footnote{The matrix $B_g$ is constructed such that $\{i,j\}$-th entry is $1$ only if the $j$-th element of vector $p$ is connected to the $i$-th node of the network, else it is $0$. A similar process is followed for the other two matrices $B_r$, $B_\ell$.} The vector of line flows at time $t$ can be expressed as $F[t] = {\Lambda} P[t]$ where $\Lambda \in \Rb^{N_e\times N_s}$ is the PTDF matrix and 
\begin{equation}\label{eq:power_injection}
P[t] := B_g p[t] + B_r w[t] - B_\ell \ell[t],
\end{equation}
denotes the vector of power injections at the nodes.
We refer the reader to \cite{modarresi2018scenario,chatzi2018lecture,cheng2005ptdf} for analytical and \cite{barbulescu2009software} for numerical approaches to compute PTDF matrices.

For ease of exposition, we equivalently represent \eqref{eq:LAED} as
\begin{subequations}
\label{eq:LAED_comp}
\begin{align}
\min_{x \in \Rb^{n_x}} \quad & c^\top x \\
\st \quad & Ax \leq b, \label{eq:laed_det} \\
& Dx + E\omega \leq f \label{eq:laed_stoch}, 
\end{align}
\end{subequations}
where $x \in \Rb^{n_x}$ is a compact representation of the decision variables $\{p_{i}[t]\}_{i \in \GG, t \in \TT}$ with $n_x = TN_g$, $\omega \in \Rb^{n_w}$ denotes the stochastic power generation by RESs $\{w_{j}[t]\}_{j \in \RR, t \in \TT}$ with $n_\omega = TN_r$, and $c, A, b, D, E$, $f$ are vectors and matrices of appropriate dimensions. In particular, we denote the dimension of $f$ by $K$, i.e., $f \in \Rb^K$. Although $\omega$ is stochastic in nature, we retain the representation as introduced in \eqref{eq:LAED} for readability. Furthermore, let $d_k^\top$ and $e_k^\top$ denote the $k$th row of the matrices $D$ and $E$, respectively. Then, the constraint \eqref{eq:laed_stoch} is equivalent to the scalar constraint 
\begin{equation}\label{eq:laed_stoch_eq}
Z(x,\omega) := \max_{k \in K} \quad d^\top_k x + e^\top_k \omega - f_k \leq 0.
\end{equation}

Our goal is to solve the above optimization problem where the uncertain constraint \eqref{eq:laed_stoch} (a compact representation of \eqref{eq:powerbalance_stochastic} and \eqref{eq:lineflow_constraint}) is modeled as a chance-constraint or via a suitable risk measure. In the former, the chance-constraint is stated as
\begin{equation}\label{eq:laed_cc}
\Pb(Z(x,\omega)\leq0) \geq 1-\alpha, \quad \alpha \in (0,1),
\end{equation}
where $\mathbb{P}$ denotes the distribution of the random variable $\omega$ and {$\alpha$ is the desired violation probability}. The chance-constraint \eqref{eq:laed_cc} ensures that {all the uncertain operational constraints will be satisfied simultaneously} with a probability of at least $1-\alpha$. 

As discussed in the Introduction, we also consider a well-established convex risk measure {\it conditional value-at-risk} (CVaR) for the uncertain constraints, which is defined as
\begin{equation}\label{eq:laed_cvar}
\CVaR{\Pb}_\alpha(Z(x,\omega)) \!:= \!\inf_{t \in \Rb} \left[\!\frac{1}{\alpha} \Eb_{\Pb}(Z(x,\omega) - t)_{+} \!+ \!t \right] \leq 0.
\end{equation}
Note that the chance-constrained program where \eqref{eq:laed_stoch} is replaced by \eqref{eq:laed_cc}, is non-convex except for a restrictive class of distributions. On the other hand, the CVaR-constrained optimization problem (with \eqref{eq:laed_stoch} replaced by \eqref{eq:laed_cvar}) is a convex conservative approximation of the chance-constrained counterpart \cite{nemirovski2006convex}. In the definition \eqref{eq:laed_cvar}, $t$ is interpreted as the smallest value such that $\Pb(Z(x,\omega) \geq t) \leq \alpha$, and $\CVaR{\Pb}_\alpha(Z(x,\omega))$ denotes the expected value of $Z(x,\omega)$ subject to $Z(x,\omega)$ exceeding $t$. Hence, CVaR captures the mean of the magnitude of the violation of the chance-constraint. 

\begin{remark}\label{remark:powerbalance}
{We formulate the power balance via the inequality constraint \eqref{eq:powerbalance_stochastic}, by following the typically accepted convention in the security-constrained OPF (SCOPF) literature (see for example \cite{capitanescu2011state}). While a number of mechanisms are available in the grid to counteract disturbances in real-time (e.g., frequency control mechanisms), their range of action is often limited and their activation may be expensive. As a result, the operator may need to take strategic actions (e.g., ramping-up traditional generation) several hours ahead of the real-time operation to ensure a greater degree of system controllability.
An alternative approach would be to model the dispatch decisions in terms of (affine disturbance) feedback policies determined by future renewable energy generation. From a methodological perspective, our formulations (discussed below) can be applied to optimize the coefficients of an affine disturbance feedback policy. In practical terms, this implies co-design of the day-ahead schedule and real-time mechanisms. While this approach is interesting, it departs from the modular architecture that is currently adopted by most operators, where scheduling and real-time operations are only coordinated via the procurement of reserves. Further investigations along these lines remains an interesting avenue for future research.} \oprocend
\end{remark}
\section{{Wasserstein Distributionally Robust LAED}}
\label{section:solution_mothodology}

We now describe the data-driven distributionally robust techniques to solve the LAED problem formulated above with chance or CVaR-constraints. Throughout, we assume that the decision-maker has access to a set of samples $\widehat{\Omega}_N := \{\widehat{\omega}_1,\widehat{\omega}_2,\ldots,\widehat{\omega}_N\}$ of the uncertain parameters with $\data_k \in \Rb^{n_\omega}$. In the LAED problem, each $\data_k$ denotes a (non-negative) vector of power generation by the RESs over an interval of length $T$. We use $[N]$ and $[K]$ to denote the sets $\{1,\ldots,N\}$ and $\{1,\ldots,K\}$, respectively.

\subsection{Distributionally Robust CVaR-Constrained LAED}
\label{ssec:drcvp}

We first consider the distributionally robust CVaR-constrainted program (DRCVP) for the LAED problem. In particular, we require the CVaR-constraint \eqref{eq:laed_cvar} to hold for a family of distributions, referred to as an {\it ambiguity set}, defined directly from observed samples via the Wasserstein metric. The corresponding optimization problem is given by%
\begin{subequations}
\label{eq:LAED_comp_drrcp}
\begin{align}
\min_{x \in \Rb^{n_x}} \quad & c^\top x \\
\st \quad & Ax \leq b, \label{eq:laed_det_drrcp} \\
& \sup_{\Pb \in \MM^\theta_N} \inf_{t \in \Rb} \left[\frac{1}{\alpha} \Eb_{\Pb}(Z(x,\omega) + t)_{+} - t \right] \leq 0 \label{eq:laed_drrcp1}, 
\end{align}
\end{subequations}
where $\MM^\theta_N$ is the Wasserstein ambiguity set defined using the samples $\widehat{\Omega}_N$. 
Specifically,
\begin{equation}\label{eq:wasserstein-set}
\MM^\theta_N := \setdef{ \mu \in \PP_1(\Omega)}{ W_1(\mu,\Pbhat_N) \leq
\theta},
\end{equation}
contains all distributions with a finite first-moment and support $\Omega$ (represented by the set $\PP_1(\Omega)$) within a distance $\theta$, measured by the Wasserstein metric, from the empirical distribution constructed from the observed samples $\Pbhat_N := \frac{1}{N}\sum^N_{i=1} \delta_{\data_i}$ ($\delta_{\data_i}$ is the unit point mass at $\data_i$). The Wasserstein metric $W_1$ is formally defined in \cite{hota2019data}. The optimization problem \eqref{eq:LAED_comp_drrcp} is potentially infinite-dimensional due to the supremum over a set of probability distributions. Next, we present a tractable finite-dimensional convex reformulation of \eqref{eq:LAED_comp_drrcp} when the support of the uncertain parameters is a {polyhedral} subset of $\Rb^{n_\omega}$. 

\begin{proposition}[Tractable DRCVP]\label{prop:drcvp}
Let the set $\Omega$ be defined as $\Omega := \{ \omega \in \Rb^{n_\omega}| G\omega \leq h \}$. Then, \eqref{eq:LAED_comp_drrcp} is equivalent to the program%
\begin{subequations}
\label{eq:LAED_drrcp_tract}
\begin{align}
\min_{x, \lambda, t, s, \eta} \quad & c^\top x \label{eq:LAED_drrcp_tract_cost}
\\ \st \quad & Ax \leq b, \label{eq:LAED_drrcp_tract_det}
\\ & \lambda \theta + \frac{1}{N} \sum^N_{i=1} s_i \le t \alpha, \label{eq:LAED_drrcp_tract_sum}
\\ & d^\top_k x \!- \!f_k \!+ t + \!(e_k - G^\top \eta_{ik})^\top \data_i + \eta_{ik}^\top h  \le s_i, \label{eq:LAED_drrcp_tract_big}
\\ & \norm{e_k - G^\top \eta_{ik}} \le \lambda, \quad \eta_{ik} \ge 0, \label{eq:LAED_drrcp_tract_eta}
\\ & t \in \Rb, \quad \lambda \geq 0, \quad s_i \geq 0, \label{eq:LAED_drrcp_tract_generic}
\end{align}
\end{subequations}
where the inequalities involving $s_i$ and $\eta_{ik}$ hold for every $i \in [N]$, $k \in [K]$, and $t$ has an analogous interpretation as in \eqref{eq:laed_cvar}. 
\end{proposition}

We present the proof in Appendix \ref{section:appendix_drcvp}. {If the support of the uncertain parameters is not known or is unbounded, i.e., $\Omega = \Rb^{n_\omega}$, then the tractable reformulation of \eqref{eq:LAED_comp_drrcp} is obtained by setting $G = 0$, $h = 0$, and without considering the decision variables $\eta$ in \eqref{eq:LAED_drrcp_tract}. Further, note that any feasible dispatch solution to problem \eqref{eq:LAED_drrcp_tract} satisfies the CVaR-constraints for all distributions within a distance $\theta$ of the empirical distribution {\it and} having a support specified by the polytope. Thus, in practical terms, if the support of the uncertainty is polyhedral and known to the decision-maker, then the optimal solution obtained by solving the problem that makes use of the information regarding the support will be less conservative than the solution obtained by solving the problem where the support is set to be unbounded.}

\subsection{Scalable Approximation of Distributionally Robust Chance-Constrained LAED via Robust Optimization}
\label{subsec:DRCCP-scale}
We recall from the earlier discussion that the CVaR-constraint \eqref{eq:laed_cvar} acts as a convex conservative approximation to the chance-constraint \eqref{eq:laed_cc}. Note, however, that the size of the above optimization problem increases with the number of samples, leading to a large computational burden. Furthermore, chance-constrained programs, and hence distributionally robust chance-constrained programs (DRCCPs) are in general non-convex. Therefore, we now present a scalable approach to approximately solve the DRCCP counterpart of the LAED problem over the ambiguity set $\MM^\theta_N$. 

The DRCCP corresponding to the problem \eqref{eq:LAED_comp} for the ambiguity set $\MM^\theta_N$ (defined in \eqref{eq:wasserstein-set}) can be stated as%
\begin{subequations}
\label{eq:LAED_comp_drccp}
\begin{align}
\min_{x \in \Rb^{n_x}} \quad & c^\top x \\
\st \quad & Ax \leq b, \label{eq:laed_det_drccp} \\
& \inf_{\Pb \in \MM^\theta_N} \Pb[Dx+E\omega \leq f] \geq 1-\alpha \label{eq:laed_drccp},
\end{align}
\end{subequations}
i.e., we require the uncertain constraints to hold jointly for all distributions in the ambiguity set.

Our approach extends an analogous approach developed in \cite{margellos2014road} for chance-constrained programs to DRCCPs. First, for each component $j$ of the uncertain parameter $\omega$, we obtain upper and lower bounds such that $\omega_j$ lies within those bounds with a high probability for all distributions in the ambiguity set. Thus, we first solve $n_\omega$ DRCCPs each with a two-dimensional decision variable. Once the bounds are computed, we construct a hyper-rectangle $\Omega^\star$ and formulate a robust optimization problem where we require the uncertain constraints to hold for all $\omega \in \Omega^\star$. The size of this robust program does not increase with the number of data points. Furthermore, it can be shown that any feasible solution of the robust optimization problem is feasible for the DRCCP \eqref{eq:LAED_comp_drccp}, i.e., the robust optimization problem is an inner approximation of the DRCCP, {by resorting to Bonferroni's inequality}.

\subsubsection{Distributionally robust bounds on each component of $\omega$} Let $w_j$ be the $j$-th component of the uncertain random vector $\omega$. Consider the DRCCP problem
\begin{subequations}
\label{eq:robust_twod}
\begin{align}
\min_{y:=(\underline{y},\overline{y})\in \Rb^{2}} \quad & \overline{y}-\underline{y} \\
\st \quad & 0 \leq \underline{y} \leq \overline{y}, \\
& \sup_{\Pb \in \MM^\theta_N} \Pb \left[ \, \omega \not \in [\underline{y}, \overline{y}] \, \right] \leq \frac{\alpha}{n_\omega}, \label{eq:robust_twod_c}
\end{align}
\end{subequations}
where $\omega$ stands for the random variable $w_{j}$ with support $\mathbb{R}_{\geq 0}$. An optimal solution $y^\star_{j}$ is such that 
\[
\Pb\left[w_{j} \in [\underline{y}^\star_{j},\overline{y}^\star_{j}] \right] \geq 1-\frac{\alpha}{n_\omega}, \quad \forall \,\Pb \in \MM^\theta_N.
\] 
However, the problem \eqref{eq:robust_twod} involves optimization over probability distributions and is infinite-dimensional. In the following, we present a finite-dimensional reformulation of \eqref{eq:robust_twod}. 

\begin{proposition}[Distributionally robust bounds]\label{prop:drccp_rob}
The optimization problem \eqref{eq:robust_twod} for the $j$-th component of $\omega$ can be equivalently stated as%
\begin{subequations}
\label{eq:robust_twod_reform}
\begin{align}
\min_{\underline{y},\overline{y},\lambda,s} \quad & \overline{y}-\underline{y}
\\ \st \quad & 0 \leq \underline{y} \leq \overline{y}, 
\\ & \lambda \theta + \frac{1}{N} \sum^N_{i=1} s_i \le \frac{\alpha}{n_\omega},\label{eq:robust_twod_reform_A}
\\ & s_i \geq 1 - \lambda \max\{0, \overline{y} - \data_{ij}\} 
\\ & s_i \geq 1 - \lambda \max\{0, \data_{ij} - \underline{y}\}  \quad \text{if $\underline{y} > 0$} 
\\ & \lambda \geq 0, \quad s_i \geq 0, \quad \forall i \in [N],\label{eq:robust_twod_reform_N}
\end{align}
\end{subequations}
where $\data_{ij}$ is the $j$-th component of the sample $\data_i$.
\end{proposition} 

{The proof of Proposition \ref{prop:drccp_rob} (which we present in Appendix~\ref{section:appendix_drccp}) relies on the exact reformulation of DRCCPs under Wasserstein ambiguity sets as stated in \cite{hota2019data}. One of the key reasons behind the intractability of this class of problems is the necessity to compute the minimum distance of the observed sample to the complement of the feasibility set (i.e., the terms comprising the summation term in \eqref{eq:indicator-zero-duality} in Appendix \ref{section:appendix_drccp}). Even for ``well-behaved" (e.g., convex) feasibility sets, the complement is usually non-convex and consequently, computing the projection to a non-convex set is often intractable. Our proof exploits the special structure of box constraints in \eqref{eq:robust_twod_c} to obtain the reformulation in \eqref{eq:robust_twod_reform}. While other classes of uncertainty sets (such as polyhedral or ellipsoidal) may lead to less conservative solutions compared to the hyper-rectangle based uncertainty sets considered here, obtaining finite-dimensional tractable reformulations for such sets is a challenging problem and remains a promising direction for future research.}

\begin{remark}
Although \eqref{eq:robust_twod_reform} is an exact reformulation of the DRCCP \eqref{eq:robust_twod}, it is still non-convex. However, since the decision variable is two-dimensional, it can be solved by nonlinear optimization solvers or via suitably designed heuristics based on line search methods up to a reasonable degree of accuracy. For the purpose of simulations, we solve the problem by adaptively updating the upper and the lower bounds $\overline y$ and $\underline y$, such that \eqref{eq:robust_twod_reform_A}--\eqref{eq:robust_twod_reform_N} is feasible at each step. \oprocend
\end{remark}

\subsubsection{Robust optimization formulation}
{Let $[\underline{y}^\star_{j}, \, \, \overline{y}^\star_{j} ]^\top$ be the optimal bounds obtained by solving \eqref{eq:robust_twod_reform} for $w_{j}$, and let $\underline{y}^\star$ and $\overline{y}^\star$ be the vectors that collect all these distributionally robust bounds. We can now solve a robust optimization problem where the constraints are required to hold for every possible realization of the uncertain vector in the hyper-rectangle 
\[
\Omega^\star := \Pi^{n_\omega}_{j =1} [\underline{y}^\star_{j},\,\, \overline{y}^\star_{j}].
\]
The resulting optimization problem is stated below.

\begin{proposition}[Scalable approximated-DRCCP]\label{prop:drccp}
The robust optimization problem
\begin{subequations}
\label{eq:LAED_rob_drccp}
\begin{align}
\min_{x \in \Rb^{n_x}} \quad & c^\top x \\
\st \quad & Ax \leq b, \label{eq:laed_det_drccp2} \\
& d_k^\top x \!+\! \max\{0, e_k^\top\!\}\, \overline{y}^\star \!+ 
\min\{0, e_k^\top\!\}\, \underline{y}^\star \leq f, \forall k \!\in\! [K] \label{eq:LAED_rob_drccp_unc},
\end{align}
\end{subequations}
where the $\max$ and $\min$ are intended as element-wise operators, is a conservative approximation of the DRCCP \eqref{eq:LAED_comp_drccp}.
\end{proposition}
}

{Note that the size of the robust optimization problem \eqref{eq:LAED_rob_drccp} is independent of the number of samples used to compute the bounds in \eqref{eq:robust_twod_reform}.} While the size of \eqref{eq:robust_twod_reform} increases with the number of samples, it is a much smaller problem (see Section~\ref{ssec:scalability}).
Through an argument analogous to \cite[Proposition 1]{margellos2014road}, it can be shown that any {\it feasible} solution of  \eqref{eq:LAED_rob_drccp}, where $[\underline{y}^\star_{j}, \,\, \overline{y}^\star_{j} ]^\top$ is {\it feasible} to \eqref{eq:robust_twod_reform}, is feasible for the DRCCP \eqref{eq:LAED_comp_drccp}. The proof is omitted in the interest of space.

\subsection{Discussion on the proposed formulations}
\label{ssec:discussion}

\subsubsection{Generality} {The tractable reformulations presented above are applicable for any Wasserstein distributionally robust chance and CVaR-constrained programs, as long as the constraint function is affine in both the decision variables as well as the uncertainty. In particular, the presented methods are applicable in formulations of the LAED problem that include curtailment factors or when the the dispatch decisions are modeled as outputs of affine disturbance feedback policies that are functions of future renewable energy generation.}

\subsubsection{{Correlations and handling joint chance-constraints}} \label{section:joint_comp}
By computing separate bounds for each component of the uncertainty, the proposed DRCCP approach ignores the correlation between them and does not exploit this information available in the data. This is the price we pay in order to obtain a scalable formulation. A related work \cite{duan2018distributionally} considered a similar robust approach to solve the DRCCP version of the OPF problem, where the upper and lower bounds were computed after normalizing the random vector by the covariance matrix, thereby preserving the spatial and temporal correlation. However, the authors in \cite{duan2018distributionally} consider individual chance-constraints {(equivalent to transforming the joint chance-constraints in \eqref{eq:laed_drccp} to individual chance-constraints by applying Bonferroni's inequality)}. As a result, the uncertainty takes the form of a scalar random variable on each constraint, which is key in ensuring the scalability of their approach. This approach from \cite{duan2018distributionally} is not applicable for joint chance constraints as it would require enumerating the vertices of a higher-dimensional hyper-rectangle to find the extremal realization of the uncertainty, which is computationally prohibitive.\footnote{We will likely encounter a similar technical challenge if we define the decision variables to be the coefficients of a suitably defined policy as opposed to the conventional power generation vector and enforce power balance constraints at all time steps.} Developing scalable approximations that preserve the correlation among the components of a random vector in joint chance-constraints is a challenging open problem and a promising avenue for future research.

{We also note a subtle difference in the way joint chance-constraints are handled in \cite{duan2018distributionally} and in our work. The application of Bonferroni's inequality in \cite{duan2018distributionally} requires the individual chance-constraints to hold with probability $\alpha/n$ where $n$ is the number of operational constraints in \eqref{eq:laed_drccp}. This scales with the size of the network, as operational constraints include line flow limits, among others. In contrast, our construction of the distributionally robust uncertainty set requires us to divide $\alpha$ by $n_\omega$ (in \eqref{eq:robust_twod}) which is the length of the random vector (number of renewable energy generators times the time horizon). Furthermore, in \cite{duan2018distributionally}, the authors construct an uncertainty set for each individual chance-constraint separately; in contrast, $\Omega^*$ is agnostic to the constraints.}

\subsubsection{Significance of the radius of the ambiguity set}
{The Wasserstein radius $\theta$ is a principled way of evaluating the trade-off between robustness and performance for the above formulations. In the robust formulation, tuning a number of different bounds on the uncertain parameters directly can be cumbersome, while the Wasserstein radius is a scalar parameter that controls the degree of distributional robustness. 

For larger values of $\theta$, we require the chance or CVaR-constraints to hold for a larger set of distributions. This is useful in instances where the number of available samples is small. While other sample-based approaches would suffer from over-fitting of the solution to the available data, under the proposed approaches, the decision-maker may choose a larger value of $\theta$ to improve the robustness of the solution to yet-to-be-realized uncertain parameters. On the other hand, when a large volume of past data is available, the empirical distribution tends to be a good representation of the uncertain parameters and the decision-maker may reduce $\theta$ to incur a smaller optimal cost. The above characteristics are highlighted in the numerical results in Section~\ref{section:numericalexperiments}.}

%%%%%%%%%%%%%%%%%%%%%%%%%%%
\subsubsection{Finite sample guarantees}\label{sec:fsg}
When the samples are drawn i.i.d. from an underlying data-generating distribution, prior works have established rigorous bounds on the Wasserstein distance between the empirical distribution and the data-generating distribution \cite{fournier2015rate,weed2019sharp}. However, the renewable energy data is not necessarily being drawn in an independent manner from any underlying distribution, as discussed earlier in the Introduction. Therefore, these guarantees are not necessarily applicable for the LAED problem considered in this work. {Nevertheless, in order to give a complete picture of the methods proposed in this work, we briefly review some of the guarantees and assumptions reported in the literature below. }

When the data-generating distribution $\Pb$ is light-tailed {(i.e., there exists $a >1$ such that $A:=\Eb_\Pb [e^{\norm{\xi}^a}] < \infty$)}, then \cite{fournier2015rate} establishes that for a given radius $\theta$, the probability with which the ambiguity set $\MM^\theta_N$ contains the data-generating distribution grows exponentially towards unity with the sample size $N$. {This result prescribes a way to select the radius $\theta$ of the ambiguity set $\MM^\theta_N$. Specifically, \cite{esfahani2018data} showed that for $n_\omega > 2$, a given $\beta \in (0,1)$, and $N$, the radius $\theta$ can be chosen as
\begin{align}\label{eq:eps-beta}
\theta_N(\beta) :=\begin{cases} \Bigl( \frac{\log(c_1 \beta^{-1})}{ c_2 N} \Bigr)^{1/n_{\omega}}, & \text{if } N \ge \frac{\log(c_1\beta^{-1})}{c_2},
\\
\Bigl( \frac{\log(c_1 \beta^{-1})}{c_2 N} \Bigr)^{1/a}, & \text{if
} N < \frac{\log(c_1 \beta^{-1})}{c_2},
\end{cases}
\end{align}
which comes with the guarantee that
\begin{align}\label{eq:fs-beta}
\Pb^N \big(W_1(\Pb,\Pbhat_N) \le \theta_N(\beta)\big) \ge 1-\beta.
\end{align}
Here, $c_1, c_2$ are positive constants that only depend on $a, A$, and $n_\omega$.} The results show that when $N$ is small, the rate of convergence is of the order of $N^{-1/a}$ where $a$ is a constant that depends on the distribution but is independent of the dimension of the uncertainty $n_\omega$. When the sample size is large, a slower convergence rate of the order $N^{-1/n_\omega}$ is observed. In a recent work \cite{weed2019sharp}, authors establish similar rates of convergence while clearly specifying, for a wide of class of distributions, the constants defining the rates of convergence. 

{The above guarantees are on the Wasserstein distance between the true and empirical distribution and as such do not depend on the distributionally robust optimization problem. Due to this generality, these theoretical guarantees are often not tight compared, for instance, to the guarantees on the probability with which the optimal solution of the scenario program satisfies original chance-constraint \cite{campi2008exact}.}\footnote{{The guarantees for the scenario program hold when the samples are drawn i.i.d. from the data-generating distribution, but do not require the latter to be a light-tailed distribution nor do they depend on any exogenous constants.}} Consequently, a much smaller value of $\theta$ (compared to the one stated in \eqref{eq:eps-beta}) often have the desired empirical performance. Section~\ref{section:numericalexperiments} provides a detailed empirical study on the out-of-sample constraint satisfaction for varying values of $\theta$, and highlight the robustness-performance trade-off. 

\section{Numerical Experiments}
\label{section:numericalexperiments}

{
\subsection{Test case description}
\label{subsection:model}

We consider the IEEE $39$-bus ``New England'' transmission grid test case \cite{pai2012energy} for the numerical experiments in this section. Three renewable sources (solar farms) have been connected to the grid, as shown in Figure~\ref{figure:IEEEnet}. The remaining traditional generators have different marginal costs depending on their type: fossil fuel, import from the grid interconnection, nuclear, and hydro (in decreasing order of cost). As the renewable generators are located in the proximity of cheap power generators, the Transmission System Operator (TSO) will strive to maximize the power flow from these buses to the rest of the grid and to ramp down the expensive sources, when possible.  All computations are carried out in MATLAB with MATPOWER \cite{ZimmermanThomas11} and MOSEK, on a personal computer with $16$ GB of memory. 

\begin{figure}[tb]
    \centering
	\includegraphics[width=0.45\textwidth]{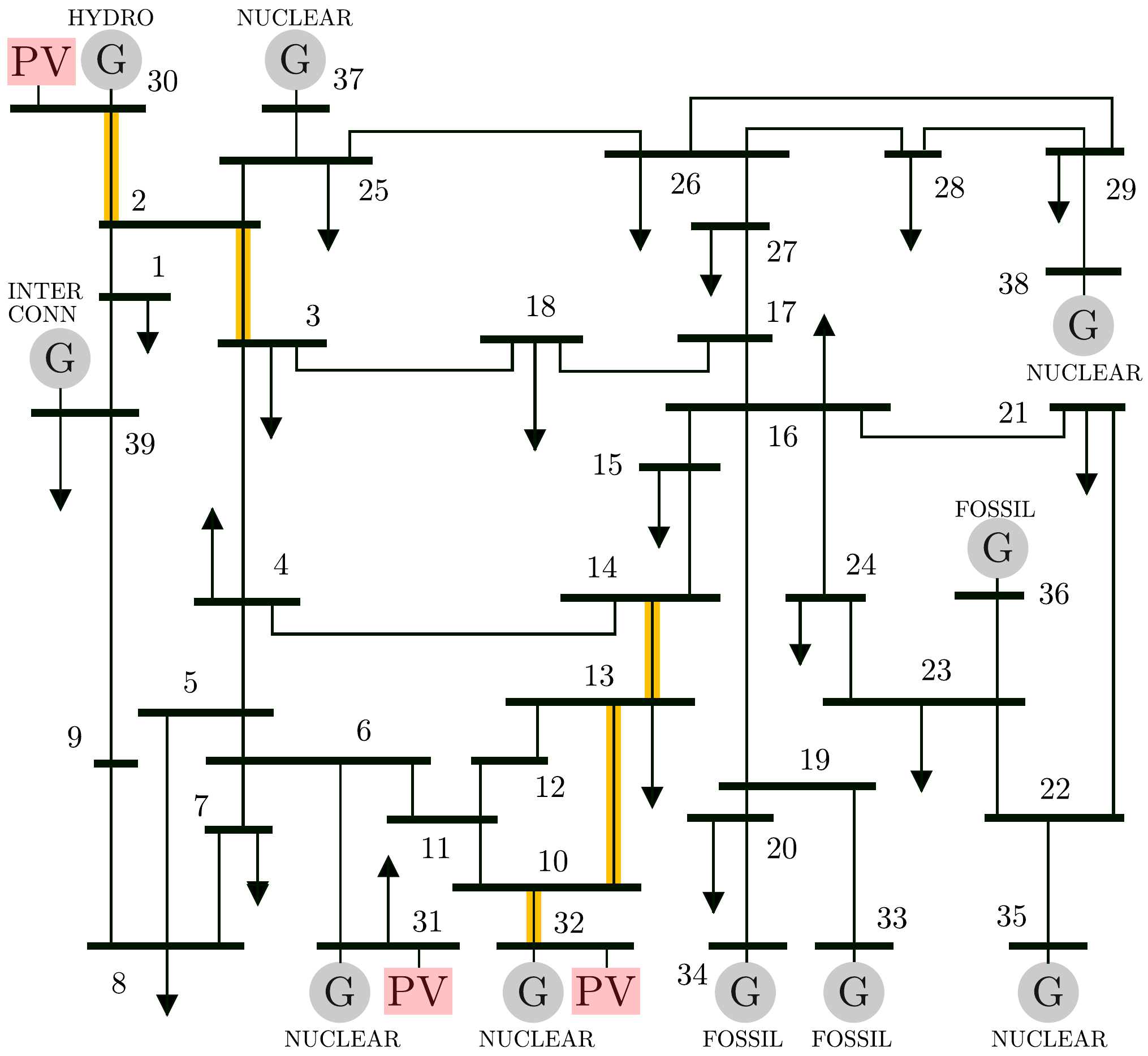}
	\caption{The benchmark IEEE $39$-bus test case, modified to include 3 PV sources, indicated in red. The power lines highlighted in yellow are prone to congestion when solar generation is abundant.}
	\label{figure:IEEEnet}
\end{figure}

\subsection{Distributionally robust dispatch with limited samples}
\label{ssec:limitedsample}

We first consider the case in which the transmission system operator has access to third party forecasts for the next-day irradiation at the locations of the solar farms. Such a forecast is often provided in the form of an \emph{ensemble} of hourly irradiation profiles. Each element of the ensemble is obtained by performing numerically intensive simulations under various meteorological models of the weather for the next day. They are, therefore, expensive to obtain and generally available in limited number. All the elements of the ensemble are to be considered as equally probable realizations and together they provide an indication of the reliability of the forecast (based on how closely they agree). 
Examples of these ensembles is reported in Figure~\ref{fig:envelope}.

For this simulation, we acquired $17$ hourly-irradiation forecasts from sources listed in Table~\ref{tab:dataset}, for the last $75$ months and for three locations in continental Europe. For the same time period and for the same locations, we also considered the satellite irradiation at hourly measurements. We then compared three possible approaches that the TSO may employ to schedule the day-ahead power generation in order to achieve a desired violation probability $\alpha$ smaller than $1$\%: 
\begin{itemize}[leftmargin=0.2cm, itemindent=0.5cm]
\setlength{\itemsep}{1pt}
\item in the \textbf{worst-case} approach, we assume that the LAED problem is solved in order to guarantee satisfaction of the grid constraints for all the possible scenarios in the ensemble;
\item in the \textbf{DRCVP} approach, the Distributionally Robust CVaR-Constrained LAED problem formulated in Proposition~\ref{prop:drcvp} is solved;
\item in the \textbf{DRCCP} approach, the Distributionally Robust Chance-Constrained LAED problem is solved employing the scalable approximation proposed in Proposition~\ref{prop:drccp}.
\end{itemize}

\begin{table}
\caption{Weather models and observations provided by meteoblue AG.}
\label{tab:dataset}
\centering
\begin{tabular}{@{}lp{45mm}@{}}
\textbf{Weather service} & \textbf{Models} \\
\midrule
meteoblue AG	& \raggedright\let\\\tabularnewline NEMS12	, NEMS12E, NEMS2-30, NEMS30, NEMS4, NMM22, NMM4 \\
NOAA (USA) & GFS05 \\
MSC (Canada) & GEM, GEM15 \\
DWD	(Germany) & ICON, ICONEU \\
Meteofrance (France) & \raggedright\let\\\tabularnewline AROME2, ARPEGE11, ARPEGE40	 \\
Met Office (UK) & UMGLOBAL	\\
KNMI (Netherlands) & HIRLAM \\
\\
\textbf{Weather service} & \textbf{Observations} \\
\midrule
METEOSAT & Satellite measurements \\ \\
\multicolumn{2}{@{}l}{Further details about these data sources are available in \cite{meteoblue}.}
\end{tabular}
\end{table}

\begin{figure}[tb]
\includegraphics[scale=1]{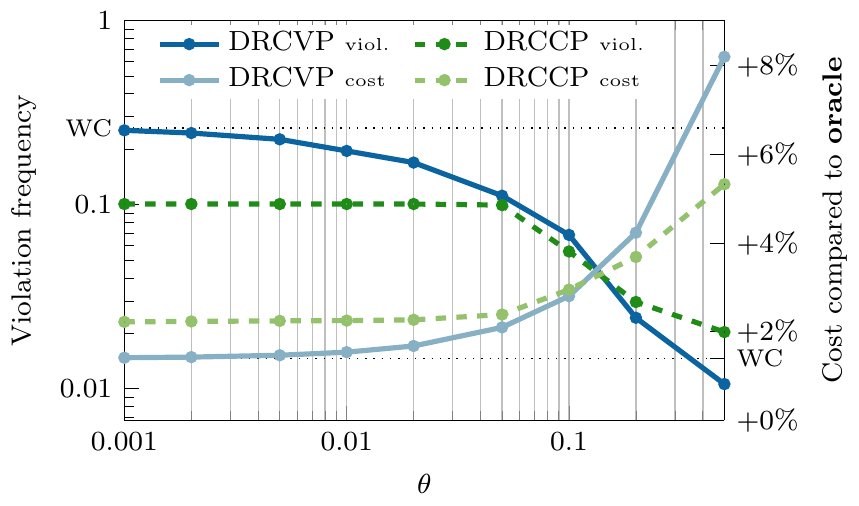}\\
\includegraphics[scale=1]{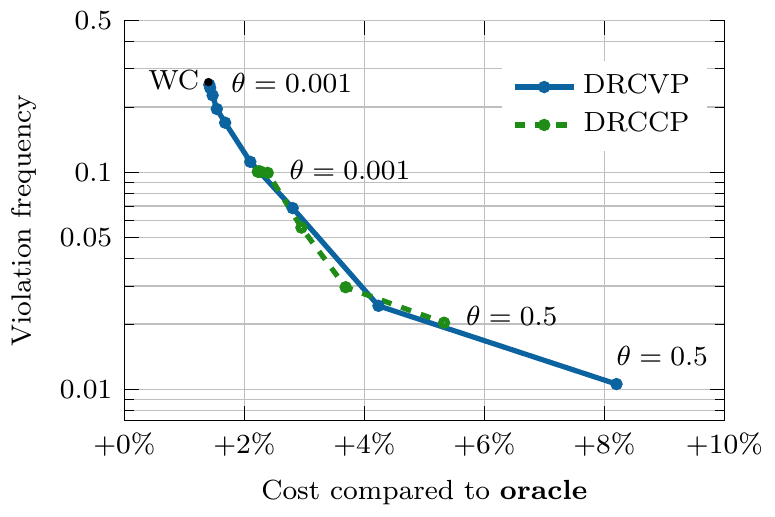}
\caption{Frequency of violation of the constraints and operational cost for the two approaches \textbf{DRCVP} and \textbf{DRCCP} for varying the radius $\theta$ of the ambiguity set. The \textbf{worst-case (WC)} approach is also marked on the axis for comparison. In all cases, the desired violation probability $\alpha$ is set to $0.01$.}
\label{fig:tradeoff}
\end{figure}

We compared these three approaches with respect to the resulting empirical frequency of constraint violations (based on the real irradiation measurements and the DC power flow solution) and with respect to the resulting operational cost for the grid. As a benchmark, we consider the \textbf{oracle} solution, the OPF that the TSO would compute if it had access to the exact irradiation profile for the next day. Figure~\ref{fig:tradeoff} shows how both DRCVP and DRCCP can be employed by the TSO to generate schedules that are safer (lower violation probability) at a small additional cost. It also shows that, given the low number of samples available in the forecast ensemble, the worst-case approach yields unsatisfactory guarantees ($26$\% violation probability). Figure~\ref{fig:tradeoff} further shows that for values of $\theta$ smaller than $0.01$, there is no significant change in the violation frequency or the optimal cost. Thus, our result shows that the solution under the DRCVP approach tends to be more robust (i.e., with a smaller violation frequency and a larger optimal cost) when $\theta$ is relatively large.

The practical implication of Figure~\ref{fig:tradeoff} is that an operator can robustify their decision to uncertain solar energy generation (as the true realized irradiation will differ from the ones included in the ensemble of forecasts) by tuning the Wasserstein radius $\theta$ which is a scalar parameter. The plot in the top panel allows a TSO to decide what value of $\theta$ should be employed in order to meet the desired violation probability and to understand the consequent cost. The plot in the bottom panel shows how both the DRCVP and the DRCCP approaches lie on the same Pareto-optimal front: none of the methods outperforms the other by producing schedules which have lower violation probability at the same cost (or, vice-versa, lower cost for the same violation probability).

\begin{figure}[tb]
	\centering
	\includegraphics[scale=1]{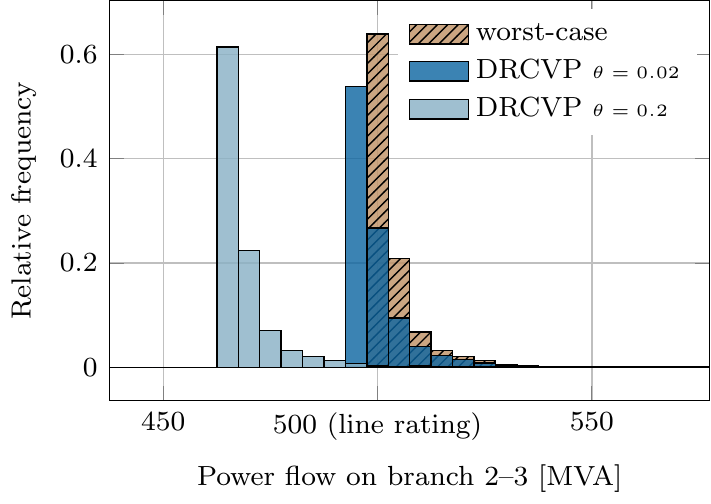}
	\caption{Statistical use of one of the most congested branches in the grid for different levels of $\theta$, in contrast to the worst-case approach.}
	\label{fig:branchuse}
\end{figure}

\begin{figure}[tb]
	\centering
	\includegraphics[scale=1]{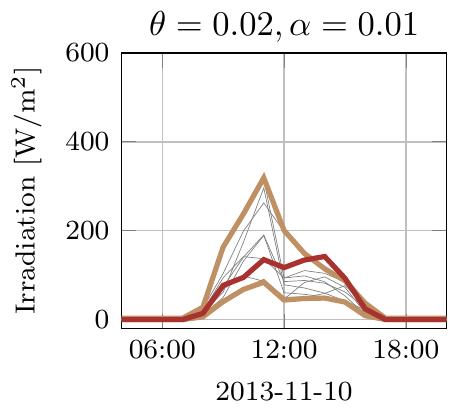}
	\includegraphics[scale=1]{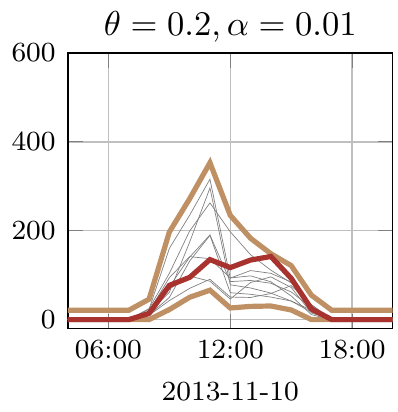} \\[2mm]
	\includegraphics[scale=1]{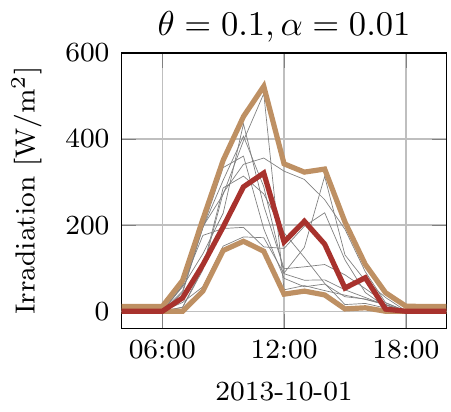}
	\includegraphics[scale=1]{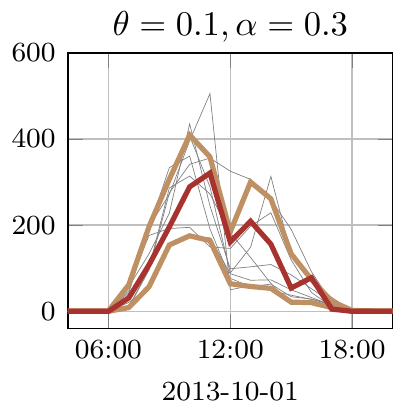}
	\caption{Illustration of the effect of the tuning parameters $\theta$ and $\alpha$. In all plots, the thin gray lines represent forecasts, the red line represents the measured irradiation, and the thick gold lines represent the lower and upper bounds $\underline{y}$, $\overline{y}$ used in the DRCCP algorithm.}
	\label{fig:envelope}
\end{figure}

A TSO may also perform a similar statistical analysis to identify the desired level of robustness (i.e., the desired $\theta$) based on the violation of specific operational constraints. Figure~\ref{fig:branchuse} provides an example of such an analysis for DRCVP (DRCCP yields similar results): the histogram shows the empirical distribution of the power flow on an critical branch for different LAED approaches. The worst-case dispatch, in which constraint satisfaction is ensured for all elements of the forecast ensemble, yields frequent violations of the line rating. In contrast, a suitably chosen value of $\theta$ produces a distribution of line flows which lies on the left of the line limit. 

An intuitive interpretation of the roles played by $\theta$ and $\alpha$ in these algorithms is offered in Figure~\ref{fig:envelope}, where we plotted the distributionally robust bounds obtained via \eqref{eq:robust_twod_reform} (for scalable approximation of DRCCP). The first row shows how a larger $\theta$ allows to be robust with respect to realizations that fall outside of the envelope defined by the few available forecasts. The second row shows how a larger $\alpha$ allows to tolerate some violation probability (in exchange for a better cost).

Thus, our results provide compelling insights on how the distributionally robust approaches can be leveraged to take dispatch decisions ahead of time and satisfy operational constraints under uncertain weather forecast (available in the form of few samples). When few predictions are available, increasing the radius $\theta$ can result in robust dispatch decisions, i.e., with a smaller likelihood of constraint violation.}
\subsection{Distributionally robust dispatch with historical data}
\label{ssec:historicalsamples}

{We now consider the case where the transmission operator has access to historical power injection data collected from the field. A TSO may be interested in using this source of data because, in contrast to third-party irradiation forecasts, they are specific of their system (for example, they factor the efficiency of their solar farms, the concurrent effect of irradiation on consumer power demand, etc.).

In order to investigate such a setting, we consider the data collected by the National Renewable Energy Laboratory (NREL) from the Sacramento Municipal Utility District (Anatolia) during the period $23$ April -- $21$ July $2012$, at a one-minute time resolution from $5$AM to $7$PM \cite{AA-LB:09}.
Hourly data have then been generated by decimating the minute-scale measurements for different intra-hour offsets, obtaining $60$ separate time series.} These historical data have two key characteristics:
\begin{itemize}[leftmargin=0.2cm, itemindent=0.5cm]
\setlength{\itemsep}{1pt}
\item In Figure~\ref{fig:corr}, we plot the correlation of solar generation between consecutive days. As illustrated, the data is highly correlated with a mean correlation coefficient of $0.99$. 
\item On a slower time-scale, Figure~\ref{fig:hist1} shows that the distribution  of solar generation at the same time of the day in two 30-day time periods are fairly different.
\end{itemize}
Based on these observations, the data cannot be assumed to be independently drawn from an underlying distribution. {As a result, we resort to empirical analysis of constraint violation as opposed relying on finite sample guarantees which hold under the assumption that samples are i.i.d., according to a true distribution. Similarly, the guarantees provided under the scenario approach \cite{modarresi2018scenario} are not necessarily applicable.} 

\begin{figure}[bt]
    \centering
    \includegraphics[scale=1]{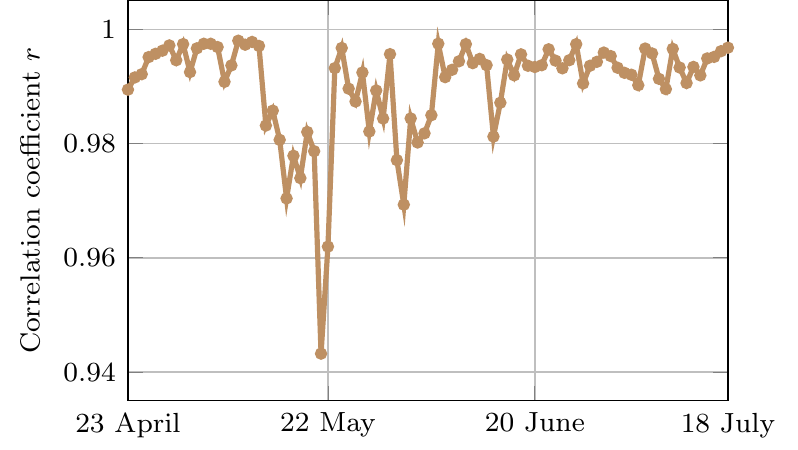}
	\caption{\footnotesize Correlation of solar generation between consecutive days.}
    \label{fig:corr}
\end{figure}

\begin{figure}[bt]
	\centering
	\includegraphics[scale=1]{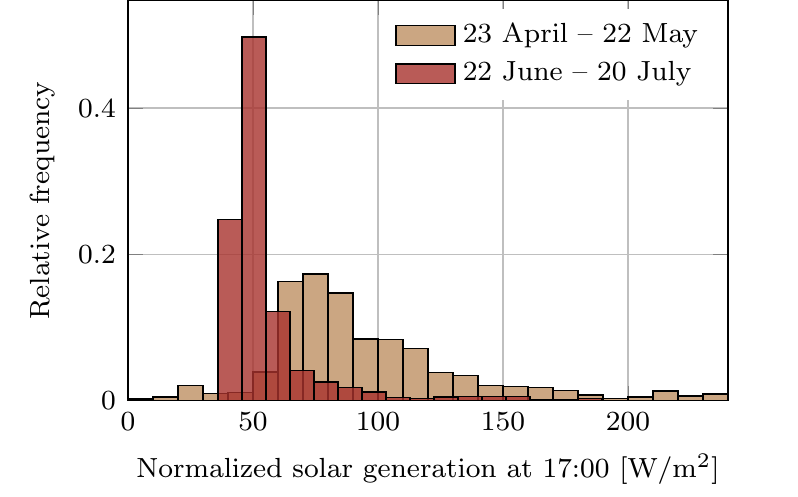}
	\caption{\footnotesize Empirical distribution of solar generation for different months.}
	\label{fig:hist1}
\end{figure}

\begin{figure}[bt]
	\includegraphics[scale=1]{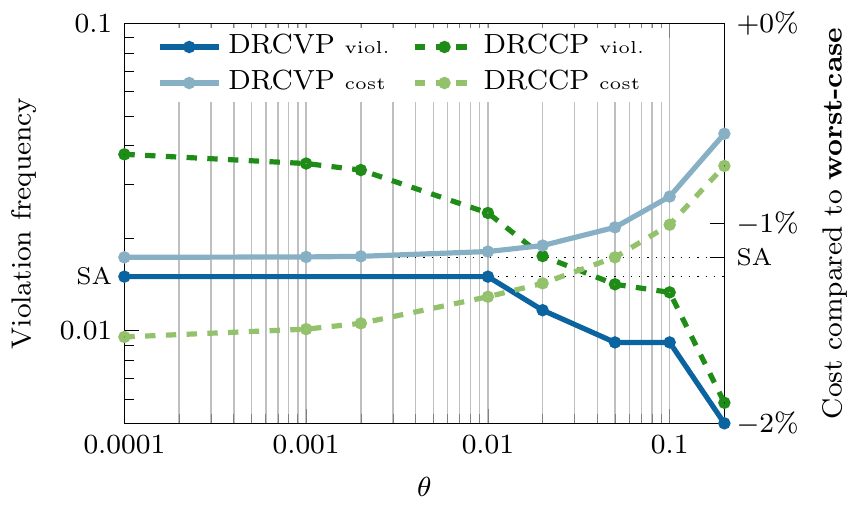}\\
	\includegraphics[scale=1]{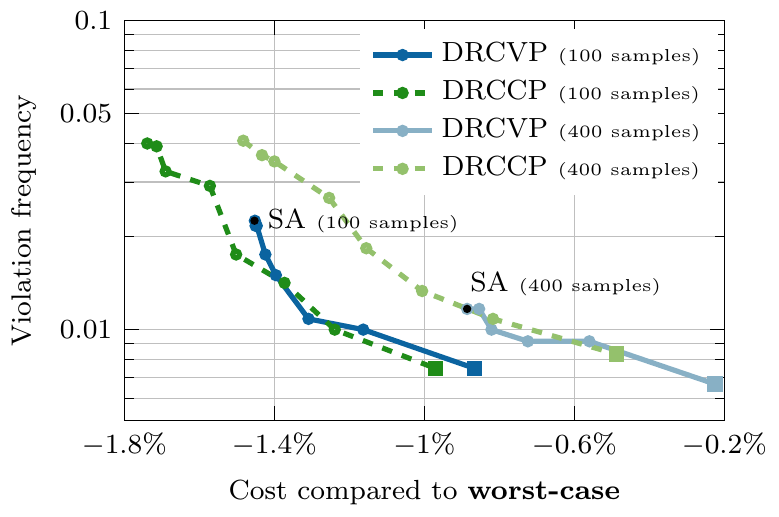}
	\caption{\footnotesize Frequency of violation of the constraints and operational cost for the two approaches \textbf{DRCVP} and \textbf{DRCCP} for $\theta$ varying between $0.0002$ and $0.2$ (the latter marked with a square), along with the Scenario approach (\textbf{SA}). In the upper plot, 200 samples are used. In all cases, $\alpha = 0.05$.}
	\label{fig:tradeoff-hist}
\end{figure}

{The LAED problem is solved using the DRCVP and DRCCP approaches, for different numbers of samples uniformly drawn from the entire dataset.
As benchmarks for these methods, we consider the \textbf{scenario} approach, where LAED is solved to ensure satisfaction of the constraints for a subset of samples of the same size, and the \textbf{worst-case} LAED solution that ensures constraint satisfaction for all $5340$ available samples. 
Note that, compared to the experiment in Section~\ref{ssec:limitedsample} where very few samples were available, the worst-case approach is now expected to be very conservative.

In order to empirically evaluate the frequency of constraint violations, a set of $1200$ samples is used as a validation dataset against which all three approaches are compared. 
Our main findings are illustrated in Figure~\ref{fig:tradeoff-hist} where we consider a violation probability of $\alpha=0.05$ and vary the parameter $\theta$ to compare the robustness (frequency of constraint violation) and performance (improvement in optimal cost compared to the worst-case solution). In the top panel of Figure~\ref{fig:tradeoff-hist}, the \% violation and improvement in optimal cost for the DRCCP and DRCVP approaches are plotted as a function of the Wasserstein radius for a training set of $200$ samples. We note that both these approaches result in an acceptable violation probability (i.e., $5$\%). The results obtained under the DRCVP are comparable to the those obtained under the scenario approach when $\theta$ is sufficiently small. Furthermore, a steep decrease in the violation frequency is observed beyond $\theta=0.01$ for both approaches.

{Nevertheless, we emphasize that our results in this section are based on the empirical frequency of violation and do not have any associated theoretical robustness guarantees since the renewable generation data is not necessarily drawn from any underlying data-generating distribution (see Section \ref{sec:fsg} for a detailed discussion).} 

The Pareto-optimal fronts for $100$ and $400$ training data samples (violation frequency is computed against the same set of $1200$ samples as before) are presented in the bottom panel of Figure~\ref{fig:tradeoff-hist}. We note that the solutions under both DRCVP and DRCCP follow similar trends and have a significant overlap for the same number of samples. As the number of samples increases, the solution becomes more robust but at a larger dispatch cost. We, however, note that increasing the number of samples is computationally expensive (especially for the scenario and the DRCVP approach -- in the next section, where we discuss the benefits of DRCCP in this regard).}

{\subsection{Scalability of the proposed approaches}
\label{ssec:scalability}

In this section, we report some observations which help us to gauge the scalability of the proposed LAED methods. To this end, we consider a modified IEEE $118$-bus network with an additional $18$ renewable sources. For this modified benchmark, the number of decision variables, constraints, runtime, and sub problems are listed in Table~\ref{tab:case3}. As before, all computations were carried out in a MATLAB environment with MOSEK on a personal computer with $16$ GB of memory.
Table~\ref{tab:case3} shows that as the number of samples increases, both the scenario approach and DRCVP scale very poorly in terms of memory footprint. While methods to remove redundant constraints exist, they typically carry a significant computational cost.
In contrast, the number of constraints of DRCCP does not increase with the number of samples.
The DRCCP approach requires the solution of a fixed number of non-convex two-dimensional subproblems (Proposition~\ref{prop:drccp_rob}), whose size increases with the number of samples. However, due to their low dimension, the total runtime remains practically constant (a few seconds), making DRCCP well suited for large problems.

\begin{table}[htb]
\centering
\caption{\footnotesize Memory footprint and computation time of different problem instances (IEEE $118$-bus network with $18$ PVs)}
\label{tab:case3}
\begin{tabular}{@{}l l P{13mm} P{13mm} P{18mm}@{}}
{} &{} & $10$ samples & $50$ samples & $200$ samples\\
\toprule
Scenario & variables & $1\,296$& $1\,296$& $1\,296$\\
& constraints & $92\,004$ &  $452\, 676$& $1\,795\,476$\\ 
& runtime &$1.00$ s &$4.77$ s & Out of memory\\
\midrule
DRCVP & variables & $1\,308$& $1\,348$& $1\,498$\\
& constraints & $92\,005$&  $452\,781$& $1\,795\,881$\\ 
& runtime &$1.02$ s &$4.96$ s & Out of memory\\
\midrule
DRCCP& variables & $1\,296$& $1\,296$ & $1\,296$\\
& constraints &  $14\,028$ & $14\,028$  & $14\,028$\\ 
& subproblems & $432$ & $432$ & $432$\\
& runtime &$2.51$ s & $2.46$ s &$2.99$ s\\
\bottomrule
\end{tabular}
\end{table}
}

\section{Conclusion}
In this paper, chance and risk-constrained multi-period economic dispatch problems are studied, and two tractable distributionally robust optimization formulations are developed in a mathematically rigorous manner. The numerical results illustrate robustness-performance trade-off of the proposed techniques. This work lays the foundation for further exploration of data-driven distributionally robust optimization techniques in power systems. 

{Several open, interesting and challenging problems have been discussed in the paper, including (1) co-design of the day-ahead schedule and real-time balancing mechanisms under uncertainty by modeling dispatch decisions as policies that depend on future renewable energy generation, (2) developing tractable reformulations for a broader class of distributionally robust uncertainty sets, and (3) developing scalable robust approximations of DRCCPs that preserve the correlations among the components of a random vector in joint chance-constraints.} Similarly, there have been limited investigations of distributionally robust semi-definite programs which are quite relevant for OPF problems. We hope this work stimulates further research in the above-mentioned topics. 
 
\appendices
\section{Proof of Proposition \ref{prop:drcvp}}
\label{section:appendix_drcvp}
\begin{IEEEproof}
We first evaluate the constraint \eqref{eq:laed_drrcp1} as
\begin{align}
& \sup_{\Pb \in \MM^\theta_N} \inf_{t \in \Rb} \left[\Eb_{\Pb}(Z(x,\omega) + t)_{+} - t\alpha \right] \nonumber
\\ = & \inf_{t \in \Rb}  \sup_{\Pb \in \MM^\theta_N} \left[\Eb_{\Pb}(Z(x,\omega) + t)_{+} - t\alpha \right] \nonumber
\\ = & \underset{t \in \Rb}{\inf} \underset{\lambda \geq 0}{\inf} [\lambda \theta^p - t\alpha \nonumber 
\\ & \textstyle \qquad +\dfrac{1}{N} \mathlarger\sum^N_{i=1} \underset{\omega \in \Omega}{\sup} [(Z(x,\omega)+t)_+ - \lambda \norm{{\omega}-\data_i}]]. \label{eq:cvar_const_full}
\end{align}
The first equality follows as a consequence of the min-max theorem in \cite{shapiro2002minimax}. The second equality is a consequence of the strong duality theorem in \cite{gao2016wasserstein}, which also shows that the infimum over $\lambda \geq 0$ is attained. On introducing auxiliary variable $s_i$ for each term in the above summation, it can be easily shown that the feasibility set of \eqref{eq:LAED_comp_drrcp} is equivalent to the set
\begin{equation}\label{eq:def_drccp_approx_reform2}
\Pi_x\left\{\begin{array}{l}  x \in \Rb^{n_x}, \\ \lambda \ge 0, \\ t \in \Rb \\ \{s_i\}_{i=1}^N \end{array} \Bigg| \begin{array}{l} Ax \leq b, \\ \lambda \theta^p + \dfrac{1}{N} \mathlarger\sum_{i=1}^N s_i \le t\alpha, \\ s_i \ge (\underset{\omega \in \Omega}{\sup} (Z(x,\omega)+t  
\\ \qquad \qquad -\lambda \norm{\omega - \data_i})_+, \forall i \in [N] \end{array} \right\},
\end{equation}
where $\Pi_x$ gives the $x$-component of the argument.

We now focus on reformulating the constraints involving $s_i \forall i \in [N]$. In particular, we have $s_i \ge 0$ and 
	\begin{align}
	s_i & \ge \underset{\omega \in \Omega}\sup  \bigl\{ \max_{k\in\until{K}} \{ d^\top_k x \! + \!e^\top_k \omega- \!f_k\}+t  -\lambda \norm{\omega -
		\data_i} \bigr\} \nonumber
	\\   & = \max_{k\in\until{K}} \bigl\{ d^\top_k x- \!f_k + t + \sup_{\omega \in \Omega} \{e^\top_k \omega - \lambda \norm{\omega - \data_i}\} \bigr\}, \nonumber
	\\ & \geq  d^\top_k x- \!f_k + t + \sup_{\omega \in \Omega} \{e^\top_k \omega - \lambda \norm{\omega - \data_i}\}, \label{eq:simplification-aff-1}
	\end{align}
	for all $k \in [K]$. In the above expressions, the second equality interchanges the $\sup$ and the $\max$. We now compute
	\begin{align}
	& \sup_{\omega \in \Omega}\{e^\top_k \omega - \lambda \norm{\omega - \data_i}\} \notag
	\\ 
	\overset{(a)}{=} & \sup_{\omega \in \Omega} \bigl\{ e^\top_k \omega -
	\sup_{\norm{z_{ik}} \le \lambda} z_{ik}^\top (\omega - \data_i) \bigr\} \notag
	\\
	\overset{(b)}{=} & \inf_{\norm{z_{ik}} \le \lambda} \bigl\{ z_{ik}^\top \data_i + \sup_{\omega \in \Omega} \{ (e_k - z_{ik})^\top \omega \} \bigr\} \notag
	\\
	\overset{(c)}{=} & \inf_{\norm{z_{ik}} \le \lambda} \bigl\{ z_{ik}^\top \data_i
	+ \inf_{\eta_{ik} \ge 0, z_{ik} = e_k  - G^\top \eta_{ik}} \eta_{ik}^\top h \bigr\} \notag
	\\
	= & \inf_{\substack{\eta_{ik} \ge 0 \\ \norm{e_k  - G^\top \eta_{ik}} \le \lambda}}
	\bigl\{ (e_k  - G^\top \eta_{ik})^\top \data_i + \eta_{ik}^\top h \bigr\}.
	\label{eq:simplification-aff-2}
	\end{align}
	Here, (a) uses the definition of the norm, (b) follows by $\inf$-$\sup$ interchange due to~\cite[Corollary 37.3.2]{rockafellar1970convex}, and (c) writes the dual form of the
	inner linear program (with 
	$\Omega = \setdef{\omega \in \Rb^{n_\omega}}{G\omega \le
		h}$). On substituting~\eqref{eq:simplification-aff-2}
	in~\eqref{eq:simplification-aff-1}, we obtain
\begin{align}\label{eq:s-ineq}
	& s_i \geq d^\top_k x- \!f_k + t +\!\!\!\!\inf_{\substack{\eta_{ik} \ge 0 \\ \norm{e_k  - G^\top \eta_{ik}} \le \lambda}} \!\!\!\!\!\bigl\{ \!(e_k - \!G^\top \!\eta_{ik})^\top\!\data_i + \eta_{ik}^\top h \bigr\}, 
\end{align}
$\forall k \in \until{K}$. It remains to be shown that the above inequality along with $s_i \ge 0$ hold if and only if there exists $\eta_{ik} \ge 0$ for all $k \in \until{K}$ such that, 
	\begin{equation}\label{eq:s-ineq-2}
	\begin{aligned}
	& s_i \geq  {d^\top_k x- \!f_k} + t + (e_k  - G^\top \eta_{ik})^\top
	\data_i + \eta_{ik}^\top h , 
	\\
	&   \norm{e_k  - G^\top \eta_{ik}} \le \lambda, \, \, s_i \ge 0, \quad \text{for all} \quad k \in \until{K}.
	\end{aligned}
	\end{equation}
	The ``if'' part in the above statement is straightforward. For the ``only if'' part consider two cases for any $k \in \until{K}$: either the $\inf$ in~\eqref{eq:s-ineq} is attained or it is not. In the former, the optimizer of the $\inf$ satisfies~\eqref{eq:s-ineq-2}. In the latter, the optimal value of $\inf$ is $-\infty$ in which case the constraint~\eqref{eq:s-ineq} is reduced to $s_i \ge 0$. Thus, one can find $\eta_{i,k}$ such that the expression on the right-hand side of the first inequality in~\eqref{eq:s-ineq-2} is negative, thereby, reducing~\eqref{eq:s-ineq-2} to $s_i \ge 0$. This concludes the proof.
\end{IEEEproof}

\vspace{-2mm}

\section{Proof of Proposition \ref{prop:drccp_rob}}
\label{section:appendix_drccp}

\begin{IEEEproof}
On drawing parallels between the notation here and that of \cite[Theorem III.1]{hota2019data}, we have $x = \begin{bmatrix} \underline{y} \quad \overline{y} \end{bmatrix}^\top, \,c = \begin{bmatrix} -1 \quad 1 \end{bmatrix}^\top,\, X = \{y \in \Rb^2 | 0 \leq \underline{y} \leq \overline{y}\},\, \xi = \omega$, and the constraint function $F(y,\omega) = \max(\omega-\overline{y},-\omega+\underline{y})$. However, we have the support of the uncertainty $\Omega = \Rb_+$ in contrast with \cite{hota2019data} where the support was $\Rb$. 

We proceed in an analogous manner as \cite{hota2019data} and evaluate 
\begin{align}
& \sup_{\Pb \in \MM^\theta_N} \Pb(F(y,\omega)> 0 ) = \sup_{\Pb \in \MM^\theta_N} \Eb_\Pb[\mathbb{1}_{\cl(\omega \in \Omega:  F(y,\omega) > 0)}] \notag
\\ & \!= \!\inf_{\lambda \geq 0} \lambda \theta +\! \frac{1}{N} \sum^N_{i=1} \! \sup_{\omega \in \Omega} [\mathbb{1}_{(\omega \in \Omega:  F(y,\omega)  > 0)}\!-\!\lambda d(\omega,\data_{ij})], \label{eq:indicator-zero-duality}
\end{align}
where $\mathbb{1}$ is the indicator function, $d(\omega,\data_{ij})$ is the Euclidean distance. The first equality follows from~\cite[Proposition 4]{gao2016wasserstein} and the second equality is a consequence of the strong duality theorem \cite[Theorem 1]{gao2016wasserstein}.\footnote{\cite[Theorem 1]{gao2016wasserstein} requires the function within the expectation to be upper semicontinuous. Since the indicator function of an open set is lower semicontinuous, we replace it with its closure. This substitution is valid due to~\cite[Proposition 4]{gao2016wasserstein}.} Now let $\Omega_1 = \cl(\omega \in \Omega: F(y,\omega) > 0)$ and $\Omega_2 = \Omega \setminus \Omega_1$. Specifically, 
\begin{equation}\label{def:omega1}
\Omega_1 :=
\begin{cases}
& [0,\underline{y}] \cup [\overline{y},\infty), \quad \text{if} \quad \underline{y} > 0, 
\\ & [\overline{y},\infty), \qquad \text{if} \quad \underline{y} = 0,
\end{cases}
\end{equation}
and thus, $\Omega_1$ is non-empty for every $y \in \Rb^2$ such that $0 \leq \underline{y} \leq \overline{y}$. For each term in the summation~\eqref{eq:indicator-zero-duality}, we introduce an auxiliary variable as
\begin{align*}
s_i & = \sup_{\omega \in \Omega} [\mathbb{1}_{(\omega \in \Omega:  F(y,\omega)  > 0)}-\lambda d(\omega,\data_{ij})]
\\ & = \textstyle \max\{\sup_{\omega \in \Omega_1} [1-\lambda d(\omega,\data_{ij})], \sup_{\omega \in \Omega_2} -\lambda d(\omega,\data_{ij}) \}.
\end{align*}

Note that if $\data_i \in \Omega_1$, the first term above is $1$ while the second term is non-positive and consequently, $s_i=1$. On the other hand, if $\data_i \in \Omega_2$, the second term is $0$. Therefore,
\begin{equation*}
s_i = \max\{0,1-\lambda \inf_{\omega \in \Omega_1} d(\omega,\data_{ij}) \}.
\end{equation*}
Finally, following the definition of $\Omega_1$ in \eqref{def:omega1}, we have
\begin{equation*}
\inf_{\omega \in \Omega_1} d(\omega,\data_{ij}) = 
\begin{cases}
& \!\!\!\!\!\max(0,\min(\overline{y} - \data_{ij}, \data_{ij} - \underline{y})), \!\!\text{ if } \!\underline{y} > 0, 
\\ & \!\!\!\!\!\max(0,\overline{y} - \data_{ij}), \qquad \text{if} \quad \underline{y} = 0.
\end{cases}
\end{equation*}
The proof follows with some rearrangement of terms.
\end{IEEEproof}

\section*{Acknowledgment}
The authors would like to thank S. Bohren at meteoblue~AG for his assistance in selecting and obtaining the most suitable meteorological data for the numerical experiments.

\bibliographystyle{IEEEtran}
\bibliography{IEEEabrv,DRCCP}

\end{document}